\documentclass{amsart}

\usepackage{amsfonts,amssymb,verbatim,amsmath,amsthm,latexsym,textcomp,amscd}
\usepackage{latexsym,amsfonts,amssymb,epsfig,verbatim}
\usepackage{amsmath,amsthm,amssymb,latexsym,graphics,textcomp}
\usepackage{graphicx}
\usepackage{color}
\usepackage{url}

\begin{document}

\newtheorem{theorem}{Theorem}[section]
\newtheorem{prop}[theorem]{Proposition}
\newtheorem{lemma}[theorem]{Lemma}
\newtheorem{cor}[theorem]{Corollary}
\newtheorem{defn}[theorem]{Definition}
\newtheorem{conj}[theorem]{Conjecture}
\newtheorem{claim}[theorem]{Claim}
\newtheorem{example}[theorem]{Example}
\newtheorem{rem}[theorem]{Remark}
\newtheorem{rmk}[theorem]{Remark}
\newtheorem{obs}[theorem]{Observation}

\newcommand{\A}{\rightarrow}
\newcommand{\C}{\mathcal C}
\newcommand\AAA{{\mathcal A}}
\newcommand\BB{{\mathcal B}}
\newcommand\DD{{\mathcal D}}
\newcommand\EE{{\mathcal E}}
\newcommand\FF{{\mathcal F}}
\newcommand\GG{{\mathcal G}}
\newcommand\HH{{\mathcal H}}
\newcommand\I{{\stackrel{\rightarrow}{i}}}
\newcommand\J{{\stackrel{\rightarrow}{j}}}
\newcommand\K{{\stackrel{\rightarrow}{k}}}
\newcommand\LL{{\mathcal L}}
\newcommand\MM{{\mathcal M}}
\newcommand\NN{{\mathbb N}}
\newcommand\OO{{\mathcal O}}
\newcommand\PP{{\mathcal P}}
\newcommand\QQ{{\mathcal Q}}
\newcommand\RR{{\mathcal R}}
\newcommand\SSS{{\mathcal S}}
\newcommand\TT{{\mathcal T}}
\newcommand\UU{{\mathcal U}}
\newcommand\VV{{\mathcal V}}
\newcommand\WW{{\mathcal W}}
\newcommand\XX{{\mathcal X}}
\newcommand\YY{{\mathcal Y}}
\newcommand\ZZ{{\mathbb Z}}
\newcommand\hhat{\widehat}
\newcommand\vfn{\stackrel{\A}{r}(t)}
\newcommand\dervf{\frac{d\stackrel{\A}{r}}{dt}}
\newcommand\der{\frac{d}{dt}}
\newcommand\vfncomp{f(t)\I+g(t)\J+h(t)\K}
\newcommand\ds{\sqrt{f^{'}(t)^2+g^{'}(t)^2+h^{'}(t)^2}dt}
\newcommand\rvec{\stackrel{\A}{r}}
\newcommand\velo{\frac{d\stackrel{\A}{r}}{dt}}
\newcommand\speed{|\velo|}
\newcommand\velpri{\rvec \,^{'}}

\title{ Packing Subgroups in Solvable Groups}

\author{Pranab Sardar}
\address{University of California, Davis}
\thanks{This research is supported by University of California, Davis.} 
\date{\today}

\maketitle

\begin{abstract}
We show that any subgroup of a (virtually) nilpotent-by-polycyclic group satisfies the
bounded packing property of Hruska-Wise(\cite{hruska-wise}). In particular, the same is true about metabelian groups
and linear solvable groups. However, we find an example of a finitely generated solvable group of derived length $3$
which admits a finitely generated subgroup without the bounded packing property. In this example the subgroup is a
metabelian retract also. Thus we obtain a negative answer to Problem $2.27$ of \cite{hruska-wise}. On the other hand,
we show that polycyclic subgroups of solvable groups satisfy the bounded packing property.
\end{abstract}

\section{Introduction}
Bounded packing was introduced by Hruska and Wise in \cite{hruska-wise}. It was motivated by the notion of
width of subgroups due to Gitik, Mitra, Rips, and Sageev \cite{GMRS} and it was already known for quasiconvex subgroups
of hyperbolic groups by the work of Mj(see Lemma 2.4, Lemma 3.3, Lemma 3.4 of \cite{mahan-CC}) although Mj did not
introduce the term `bounded packing'. However, Hruska and Wise reproved bounded packing
for quasi-convex subgroups of hyperbolic groups in \cite{hruska-wise}. Their main theorem states that the relatively
quasi-convex subgroups of relatively hyperbolic groups satisfy bounded packing under mild restrictions. 
Although bounded packing holds for all finitely generated subgroups of free groups, nilpotent groups and
polycyclic groups (\cite{yang}) it is not the case in general. For example, bounded packing for finitely
generated subgroups of hyperbolic groups is false in general as is shown by Hruska-Wise.
However, there are few examples where bounded packing is known to fail as is asserted by these authors. This
paper is motivated by the following problem that these authors ask.

\smallskip
{\bf Problem} (Question 2.13 of \cite{hruska-wise} ): {\em Let S be a (countable) solvable group. Does every subgroup of
S have bounded packing?}
\smallskip

The authors conjectured that the answer is {\em yes} for polycyclic groups. Later this was verified by W. Yang(\cite{yang})
who proved it more generally for any separable subgroup of a countable group.
However, in this paper we give a negative answer for general solvable groups by producing a finitely generated
solvable group of derived length $3$ which admits a finitely generated subgroup without the bounded packing property.
This subgroup is also a retract and it, therefore, answers Problem $2.27$ of \cite{hruska-wise} negatively. However, we show
that finitely generated subgroups of (virtually) metabelian groups have bounded packing. More generally, we show the
following:

\smallskip

{\bf Theorem \ref{nice-thm}} {\em Suppose we have an exact sequence of groups
$1 \rightarrow N \rightarrow G \stackrel{\pi}{\rightarrow} Q \rightarrow 1$ where $N$ is nilpotent and $Q$ is polycyclic.
Then any (finitely generated) subgroup of $G$ satisfies the bounded packing property.}

\smallskip

We remark that metabelian groups are not subgroup separable in general.
We mention a new set of examples of finitely generated groups admitting finitely generated subgroups without
bounded packing using the work of Yves Cornulier \cite{cornu}.  
\smallskip

\noindent {\bf Acknowledgments:} The author would like to thank Chris Hruska, and Mahan Mj for
useful discussions. Heartfelt thanks are also due to Michael Kapovich for many many motivating, and helpful discussions
over the last two years and his constant support.

\section{Bounded Packing Property}

\begin{defn}(\cite{hruska-wise})
Suppose $G$ is a countable group with a proper, left invariant metric $d$. Let $H\leq  G$ be a subgroup.
We say that $H$ has bounded packing in $G$ (with respect to $d$) if for any $D>0$ there is a number
$n=n(G,H,D)$ such that given any collection of distinct cosets $I= \{gH: g\in G\}$ of $H$ in $G$ such
that $d(g_1H, g_2H)\leq D$ for all $g_1H,g_2H \in I$, we have $|I|\leq n$.
\end{defn}

By Lemma 2.2 of \cite{hruska-wise} we know that bounded packing  for a subgroup of a countable group with a proper
left invariant metric is independent of the particular metric. By the work of Higman, Neumann and Neumann(\cite{hnn})
we know that any countable group can be embedded in a $2$-generated group; it follows that any countable 
group admits a left invariant proper metric and thus bounded packing of subgroups makes sense for
all countable groups.

The following lemma is a summary of elementary properties in connection with bounded packing which are proved
by Hruska-Wise in the section $2$ of \cite{hruska-wise}. We will make repeated use of it later.

\begin{lemma}\label{lemma1} (\cite{hruska-wise}) Let $G$ be a countable group with a proper left invariant metric.
Then the following are true.
\begin{enumerate}
\item Every finite subgroup of $G$ has bounded packing in $G$. Also every finite index subgroup of $G$ has bounded packing in $G$.
\item Every normal subgroup of $G$ has bounded packing in $G$.
\item For any sequence of subgroups $K\subset H\subset G$ if $K$ has bounded packing in $H$ and 
$H$ has bounded packing in $G$, then $K$ has bounded packing in $G$.
\item If $H,K$ are two subgroup of $G$ both of which have bounded packing in $G$ then so does their intersection
$H\cap K$.
\item Suppose we have a short exact sequence of groups $1\rightarrow N\rightarrow G\stackrel{\pi}{\rightarrow} Q\rightarrow 1$.
Let $Q_1$ be a subgroup of $Q$. Then $Q_1$ has bounded packing in $Q$ if and only if $\pi^{-1}(Q_1)$ has the bounded packing in $G$.
\end{enumerate}
\end{lemma}

\begin{theorem}\label{nilp-thm}(Theorem 2.12 of \cite{hruska-wise})
If $G$ is a virtually nilpotent group then each subgroup of $N$ has bounded packing.
\end{theorem}

\smallskip
\noindent
{\bf Bounded Packing for Intersection of Subgroups}

The following theorem was essentially proved by Yang in \cite{yang}. 
To have it in this form requires a trivial modification of the proof of Yang \cite{yang}
for his theorem on bounded packing of separable subgroups. We include a proof for the sake of
completeness. However, although the statement is interesting we don't know of any remarkable
application.

\begin{theorem}(\cite{yang}) \label{th1}
Suppose $G$ is a countable group with a proper left invariant metric.
Let $\{H_{\alpha}\}_{\alpha\in I}$ be any collection of subgroups of $G$ such that each $H_{\alpha}$ has bounded
packing in $G$. Let $K=\cap_{\alpha\in I} H_{\alpha}$. Then $K$ has bounded packing in $G$.
\end{theorem}

$Proof:$ Let $D>0$ be given and let $\{a_i K\}_{i\in \mathbb N}$ be any collection of 
left cosets of $K$ such that $d(a_iK, a_jK)\leq D$ for all $i,j\in \mathbb N$. Then all these cosets are contained in a finite
collection of double cosets $\{K, Kb_1K, Kb_2K,...,Kb_nK\}$ where $d(1,b_i)\leq D$ and $b_i\not\in K$ for all $1\leq i \leq n$,
for some $n \in \NN$.
Let $S:=\{g\in G\setminus K: d(1,g)\leq D\}$. This is a finite set since the metric of $G$ is proper. Note that
$b_i\in S$ for all $1\leq i\leq n$.

For each $t\in S$ choose $H_t:= H_{\alpha_t}$, $\alpha_t \in I$ such that $t\not \in H_{\alpha_t}$.
Let $H=\cap_{t\in S} H_t$. We know that $K\subset H$. Note that $H$ has bounded packing in $G$ by Lemma \ref{lemma1}$(4)$.

We claim that all the distinct cosets of $K$ from the above collection are contained
in distinct cosets of $H$. Otherwise, suppose $a_iK,a_jK$ are distinct cosets of $K$ and $a_iK\subset gH$, $a_jK\subset gH$
for some $g\in G$.
This implies that $g^{-1}a_iK, g^{-1}a_jK\subset H$. This in turn implies that $Ka^{-1}_ia_jK\subset H$.
However, $Ka^{-1}_ia_jK=Kb_iK$ for some $b_i$ whence $b_i\in H$. This is a contradiction since $b_i\in S$.

We thus have a collection of distinct cosets $\{a_iH\}$ of $H$ such that $d(a_iH,a_jH)\leq D$ for all $i,j \in I$.
Since $H$ has bounded packing in $G$ the set $I$ is finite. This completes the proof. $\Box$

\begin{cor}\label{yang-thm}(\cite{yang})
Separable subgroups in any countable group satisfy the bounded packing property.
\end{cor}

$Proof:$ We know that finite index subgroups have bounded packing by Lemma \ref{lemma1}$(1)$. 
Hence an intersection of any collection of finite index subgroups will have bounded packing by the above
theorem. $\Box$

\begin{cor}\label{yang-thm2}(\cite{yang})
Subgroups of polycyclic groups have bounded packing.
\end{cor}
$Proof:$ We know that any subgroup  of a polycyclic group is separable by a  result of Hirsch (\cite{hirsch}). $\Box$


\section{Bounded Packing in virtually Abelian-by-Polycyclic Groups}

In this section we prove that finitely generated subgroups of virtually abelian-by-polycyclic groups
have bounded packing (Proposition \ref{main-thm}). To prepare for the proof
we need the following set of lemmas. 

The proof of the following lemma is motivated by a similar result due to
Cornulier(Lemma 3.18 of \cite{cornu}). Similar ideas are used throughout the paper,
most prominently in Lemma \ref{mylemma2}.

\begin{rem}\label{imp-rem}
\begin{enumerate}
\item To keep the notation simple, we will always denote the naturally defined quotient maps by $\pi$ in all the proofs,
unless otherwise needed.
\item Note that all the groups that we are concerned about sit in the middle of a short exact sequence.
If instead we assume that the groups are ``virtually'' have this property then the conclusions will still hold by
Lemma \ref{lemma1} $(1)$ and $(3)$. However, we omit the adverb ``virtually" for simplicity.
\end{enumerate}
\end{rem}

\begin{lemma}\label{mylemma2}
Suppose we have an exact sequence of groups $1 \rightarrow N \rightarrow G \stackrel{\pi}{\rightarrow} Q \rightarrow 1$
where $N$ is abelian. Suppose $Q_1:= \pi(H)$ has bounded packing in $Q$. Let $K=H\cap N$ and $H_1:=H/K$.
Suppose $H$ does not have bounded packing in $G$. Then there is a split exact sequence of the form
$$ 1 \rightarrow M \rightarrow G_1 \stackrel{\pi}{\rightarrow} H_1 \rightarrow 1$$
where $M$ is a finitely generated submodule of the $\ZZ[H_1]$-module $N/K$,  and $H_1$ does not have bounded packing in $G_1$.
\end{lemma}

$Proof:$
By Lemma \ref{lemma1}(5) $\pi^{-1}(Q_1)$ has bounded packing in $G$. Hence, by Lemma \ref{lemma1}(3) $H$ does not have bounded packing
in $G^{'}:=\pi^{-1}(Q_1)$ since $H\subset G^{'}$. Now going modulo $K$ we get an exact sequence
$$1 \rightarrow N/K \rightarrow G^{'}/K \stackrel{\pi}{\rightarrow} Q_1 \rightarrow 1.$$
By Lemma \ref{lemma1}(5) $H_1=H/K$ does not have bounded packing in $G^{'}/K$.
Therefore, there are infinitely many cosets $g_iH_1$ which are pairwise close where $g_i\in N/K$.
It is clear that we can assume all the $g_i$'s to be in a finitely generated submodule, say $M$, of $N/K$ using the properness
of the metric on $G^{'}/K$. Thus the bounded packing of $H_1$ fails in the subgroup, say $G_1$, of $G^{'}/K$ generated
by $M$ and $H_1$. Now we have a natural exact sequence $1 \rightarrow M \rightarrow G_1 \stackrel{\pi}{\rightarrow} Q_1 \rightarrow 1$.
Note that $H_1$ surjects onto $Q_1$ and intersects $M$ trivially. Hence $\pi$ restricted to $H_1$ is an isomorphism onto $Q_1$.
Thus this is a split exact sequence and the action of $H_1$ (or equivalently $Q_1$) gives a finitely generated module structure
on $M$. $\Box$

\begin{prop}
Suppose a group $G$ fits into an exact sequence of groups of the form
$1\rightarrow A\rightarrow G\rightarrow Q\rightarrow 1$ where $A$ is a finitely generated abelian group and
all (finitely generated) subgroups of $Q$ satisfies the bounded packing property. Then all (finitely
generated) subgroups of $G$ satisfy the bounded packing property.
\end{prop}
$Proof:$ Let $H$ be any (finitely generated) subgroup of $G$. By the given conditions $\pi(H)$ has bounded packing
in $Q$. Let $K=H\cap N$ and $H_1=H/K$ as in the previous lemma. If $H$ does not have bounded packing in $G$ then
by the above lemma there is a split exact sequence
$1 \rightarrow M \rightarrow G_1 \stackrel{\pi}{\rightarrow} H_1 \rightarrow 1$
where $M$ is a finitely generated $\ZZ[H_1]$-module of $N/K$ and $H_1$ does not have bounded packing in $G_1$.
In particular, $M$ is a finitely generated abelian group. 

{\bf Claim:} $H_1$ is separable in $G_1$.

Note that this gives a contradiction due to Corollary \ref{yang-thm} and proves the lemma.

{\em Proof of the claim:} For all integer $k\geq 2$ $k.M$ is normal subgroup of $G_1$ and clearly
the image of $H_1$ in $G/k.M$ has finite index since $M/k.M$ is a finite group. Therefore,
the subgroup, say $H_k$, of $G_1$ generated by $ H_1$ and $k.M$ is of finite index in $G_1$. 
Therefore, it is enough to observe that $H_1=\cap_{k\geq 2} H_k$. However, this is equivalent to 
showing $\cap_{k\geq 2}k.M=(0)$ which is clear since $M$ is a finitely generated abelian group. $\Box$

\smallskip
\noindent
\begin{rem} 
We will see below (see Proposition \ref{example1}) that the assumption of finite generation of $A$ is necessary. 
\end{rem}

\begin{theorem}\label{main-thm}
Suppose we have an exact sequence of groups $1 \rightarrow N \rightarrow G \stackrel{\pi}{\rightarrow} Q \rightarrow 1$
where $N$ is abelian and $Q$ is polycyclic. Then any (finitely generated) subgroup of $G$ satisfies the bounded packing
property.
\end{theorem}

$Proof:$ Let $H\subset G$ be a subgroup. If $H\subset N$ then we are done by Lemma \ref{lemma1}(2),(3) because $N$ is abelian and
normal in $G$. So we may suppose $\pi(H)$ is nontrivial. By Corollary \ref{yang-thm2} $\pi(H)$ has bounded packing in $Q$,
since $Q$ is polycyclic. Since subgroup of a polycyclic group is again polycyclic $\pi(H)$ is polycyclic. Note that
it is isomorphic to Let $H_1=H/K$ where $K:=N\cap H$. Now, if $H$ does not have
bounded packing in $G$ then, by Lemma \ref{mylemma2}, we can find a finitely generated module $M$ over $H_1$ such that
$H_1$ does not have bounded packing in the semi-direct product $G_1:=H_1\ltimes M$ where $M$ is a finitely generated
$R= \ZZ[H_1]$-module. 

{\bf Claim:} $H_1$ is separable in $G_1$.

First note that the claim gives a contradiction by Corollary \ref{yang-thm} and proves the theorem.
For the proof of the claim we make use of the following result of Roseblade(\cite{roseblade2}) (and independently
by Jategaonkar(\cite{jateg})).

{\em * Finitely generated abelian-by-polycyclic groups are residually finite.}

{\em Proof of the claim:}  
By the result of Jategaonakar and Roseblade $G_1$ is residually finite. 
Therefore, $G_1$ has a plenty of finite index subgroups whose intersection is trivial. 
Intersecting finite index subgroups of $G_1$ with $M$ we get finite index subgroups of $M$. Since $G_1$ is finitely
generated there are only a finite number of finite index subgroups of a fixed index. Therefore, there is a sequence of 
finite index subgroups, say $\{M_k\}$ of $M$ which are normal in $G_1$. As in the proof of above proposition,
we now go modulo $M_k$ to get a split exact sequence
$$1\rightarrow M/M_k\rightarrow G_1/M_k\rightarrow H_1\rightarrow 1$$
Since $M/M_k$ is finite, $H_1$ is of finite index in $G_1/M_k$.
Thus the group generated by $M_k$ and $H_1$ is of finite index in $G_1$. Denote this group by $H_k$.
Since $G$ is residually finite it follows that $\cap M_k=(0)$. It follows then that $\cap H_k =H_1$.
This means that $H_1$ is separable in $G_1$.
$\Box$

\begin{cor}\label{th2}
Suppose a finitely generated group $G$ fits into an exact sequence of groups of the form
$1\rightarrow A\rightarrow G\stackrel{\pi}{\rightarrow} Q\rightarrow 1$ where $A$ is any abelian group and
all (finitely generated) subgroups of $Q$ satisfies the bounded packing property. Then all (finitely
generated) subgroup $H$ of $G$ such that $\pi(H)$ is a polycyclic group, satisfies the bounded packing property.

In particular, all polycyclic subgroups of $G$ satisfies bounded packing.
\end{cor}

\begin{theorem}\label{nice-thm}
Polycyclic subgroups of solvable groups satisfy bounded packing.
\end{theorem}

$Proof:$ This follows from Corollary \ref{th2} by induction on the derived length of solvable groups. $\Box$

\section{Bounded Packing in Nilpotent-by-Polycyclic Groups}

We are now ready to prove our main theorem.

\begin{theorem}
Suppose a group $G$ fits into an exact sequence of groups of the form
$1\rightarrow N\rightarrow G\stackrel{\pi}{\rightarrow} Q\rightarrow 1$ where $N$ is any nilpotent group and
$Q$ is a polycyclic group. Then all subgroups of $G$ satisfies the bounded packing property.
\end{theorem}

$Proof:$ The technique of proof for this theorem is similar to the one used for Lemma \ref{mylemma2}.
Let $H$ be a subgroup of $G$ without bounded packing. We may assume that $\pi(H)$ is nontrivial by
Theorem \ref{nilp-thm} and Lemma \ref{lemma1}(2). We shall show a contradiction from this. As in the proof of 
Lemma \ref{mylemma2} we reduce to the case $\pi(H)=Q$. Now, we shall use induction on the length of the central
series of $N$. When $N$ is abelian we are done by Theorem \ref{main-thm}. Suppose $N$ is a nonabelian nilpotent
group. Let $C$ be its center. Then $N/C$ has a strictly smaller central series than that of $N$.

Consider the short exact sequence $$1\rightarrow N/C\rightarrow G/C\stackrel{\pi}{\rightarrow} Q\rightarrow 1.$$
The image of $H$ in $G/C$ has bounded packing by induction. Call it $\bar{H}$. Hence $H$ fails to have bounded packing in
the inverse image of $\bar{H}$ in $G$. Call it $G_1$. Restricting the map $\pi$ we have an exact sequence
$$1\rightarrow K\rightarrow G_1\stackrel{\pi}{\rightarrow} Q\rightarrow 1.$$
Now, $G_1$ is generated by $C$ and $H$ and $K$ is a subgroup of $N$. Let $L=K\cap H$.

We observe that $L$ is a normal subgroup of $G_1$ because $L$ is normal in $H$ and $L$ is a subgroup of $N$ which
is normalized by $C$- the center of $N$. We are using here that
$G_1$ is generated by $H$ and $C$. Now, go modulo $L$ to get the following exact sequence.
$$1\rightarrow K/L\rightarrow G_1/L\stackrel{\pi}{\rightarrow} Q\rightarrow 1$$
Since $H$ does not have bounded packing in $G_1$, $H/L$ can not have bounded packing in $G_1/L$ by Lemma \ref{lemma1}(5).
However, note that $\frac{H}{L}\cap \frac{K}{L}$ is trivial and $H/L$ surjects onto $Q$. Thus $H/L$ is polycyclic.
This means $H/L$ has bounded packing in $G_1/L$ by Theorem \ref{nice-thm}. This is a contradiction. $\Box$

\begin{cor}
All subgroups of a linear solvable group satisfies the bounded packing property.
\end{cor}

By linear group we mean a subgroup of $GL(n,F)$ for some integer $n\in \NN$ and a field $F$. 

$Proof:$ We know that all solvable linear groups are virtually nilpotent by abelian. Hence we can
apply Remark \ref{imp-rem}(2) along with the above theorem. $\Box$

\begin{rem}
It is not true that all linear solvable groups are subgroup separable e.g. the solvable Baumslag-Solitar
group $BS(1.2)=<a,b| aba^{-1}=b^2>$ is not subgroup separable, although there are a few instances where
it is true, see e.g. \cite{alperin-farb}.
\end{rem}

\section{An Example}

\subsection{ The example on solvable groups}

\begin{prop}\label{example1} 
There is a finitely generated solvable group of derived length $3$ which admits a finitely generated
subgroup without the bounded packing property.
\end{prop}

For the proof of this proposition we need the following standard lemma.

\begin{lemma} \label {ref1}
Let $H$ be a subgroup of a group $G$. Suppose $H$ is acting on a vector space $V$ over $\mathbb Q$.
Then $G$ acts in a natural way on $\oplus_{n\in S} V$ where $S$ is a set of left
coset representatives of $H$ in $G$ such that on the copy of $V$ corresponding to $H\in S$ the $H$-action
is the given one.
\end{lemma}

{\em Proof of the proposition:}
Let $T$ be the group of positive rational numbers under multiplication. Then it acts naturally on
the group $(\mathbb Q,+)$.

Note that $T$ is isomophic to the free abelian group of countably infinite rank. The set of primes
forms a basis $B$. We can embed $T$ in $Q:=\ZZ \wr \ZZ$ as the obvious normal subgroup with infinite cyclic quotient.
Now using Lemma $\ref{ref1}$ we can obtain an action of this wreath product on $W:=\oplus_{n\in \mathbb Z} \mathbb Q$. 

Consider the semi-direct product $G:=W\rtimes Q$. It is clear that $W$ is 
a cyclic module over $\ZZ[Q]$ since $\mathbb Q$ is a cyclic module over $\ZZ[T]$.
Hence, $G$ is a finitely generated group.

{\bf Claim:} $Q$ does not have bounded packing in $G$. 

{\em Proof of the claim:} 
$T$ is a normal subgroup of $Q$. Hence
it suffices to show that $T$ does not have bounded packing in $G$. Now, $T$ is acting on $W$ and leaves
the $0$-th copy of $\mathbb Q$ invariant. We note that the cosets of $T$ labelled by $\mathbb Q$
are all close to each other, for suppose $rT$ and $sT$ are two distinct cosets where $r>s$. We need to find an
element of the double coset $T(r-s)T$ which is uniformly close to the identity. Since $t:=r-s>0$,
clearly we can choose $t.(r-s).t^{-1}\in T(r-s)T$. This gives us the element $1\in \mathbb Q$. Since this
does not depend on $r,s$ we are done. $\Box$

\begin{rem}
Note that bounded packing for a general countable metabelian group is false since $T$ does not
have bounded packing in $T\ltimes \mathbb Q$.
\end{rem}

\subsection{ Some nonsolvable examples}

Suppose $G$ is a countable group with a left invariant proper metric and $H$ is a subgroup of $G$.
The for all $g_1,g_2\in G$, $d_G(g_1H,g_2H)\leq D$ if and only the $d_G(Hg^{-1}_1g_2H, 1)\leq D$.
Thus if $H$ has only finitely many double cosets in $G$ but infinite index in $G$ then clearly $H$ fails to have bounded packing.
In other words, {\em if an infinite index subgroup of a group has finite bi-index (as defined in \cite{cornu})
then it fails to have bounded packing.}

\begin{lemma} Suppose $N$, $Q$ are two finitely generated groups and $Q$ acts $2$-transitively on an infinite set $X$.
Let $G= N \wr_X Q$ be the generalized restricted wreath product of $Q$ with $N$.
Then $Q$ does not have bounded packing in $G$.
\end{lemma}

Recall that the generalized wreath product here is
the semi-direct product for the action of $Q$ on $\oplus_{x\in X}N$ by permuting co-ordinates labelled by $X$.
When $X=Q$ and the $Q$-action is by left multiplication, we get the usual restricted wreath product $N\wr Q$.

$Proof:$ Fix $x\in N$ and denote by $x_t$, ${t\in X}$ the image of $x$ in the copy of $N$ correspoding to the co-ordinate $t\in X$.
Then clearly $x_t Q$ is an infinite collection of cosets of $Q$ which are all pair-wise close to each other. $\Box$

\begin{cor}
We have a split exact sequence of groups $1\rightarrow N\rightarrow G \rightarrow \ZZ*\ZZ/2\rightarrow 1$
where $G$ is finitely generated such that $\ZZ*\ZZ/2$ does not have bounded packing in $G$.
\end{cor}
$Proof:$ Note that $\ZZ*\ZZ/2$ acts naturally on $X:=\ZZ$ where $\ZZ$ acts by translation and $\ZZ/2$ 
exchanges $0$ and $1$ and keeps everything else fixed. This action is clearly $2$-transitive. Hence, we
can define $G$ to be the generalized wreath product $\ZZ\wr_{\ZZ} \ZZ*\ZZ/2$. Then the above lemma finishes the proof.

\begin{rem}
Several other examples of groups admitting $2$-transitive actions on infinite sets are given by Cornulier (
see Examples 3.4, 3.5, 3.6 in \cite{cornu}). 
\end{rem}

\bibliography{bdd-packing.bib}
\bibliographystyle{amsalpha}

\end{document}